%% file: CDC2018.tex
\documentclass[letterpaper, 10 pt, conference]{ieeeconf}  % Comment this line out if you need a4paper

\IEEEoverridecommandlockouts                              % This command is only needed if 
\newtheorem{theorem}{Theorem}

\newtheorem{procedure}{Procedure}

\newtheorem{example}{Example}

\usepackage{comment}

\input Maths.tex

\usepackage[cmex10]{amsmath}
\usepackage{amsfonts}
\usepackage{amssymb}
\usepackage{mathtools}
\usepackage{nicefrac} 
\usepackage{mathrsfs}
\usepackage{cleveref,graphicx}
\usepackage[toc,acronym,nonumberlist]{glossaries}
\usepackage{enumerate,color}
\usepackage{overpic}
\usepackage{subfigure}

\usepackage{pgfplots} 

\pgfplotsset{compat=newest} 
\pgfplotsset{plot coordinates/math parser=false} 
\newlength\fheight 
\newlength\fwidth 

\usepackage{tikz}
\usetikzlibrary{arrows,shapes,positioning}

\usepackage{cleveref}
\usepackage{autonum}
\crefname{problem}{problem}{problems}
\crefname{proposition}{proposition}{propositions}
\crefname{procedure}{procedure}{procedures}
\crefname{assumption}{assumption}{assumptions}

\title{\LARGE \bf
Sensitivity Function Trade-offs for Networks with a String Topology
}

\author{Richard Pates and Kaoru Yamamoto% <-this % stops a space
\thanks{R. Pates is with the Department of Automatic Control, Lund University, Box 118, {SE-221 00 Lund}, Sweden. This author is a member of the LCCC Linnaeus Center and the ELLIIT Excellence Center at Lund University.} 
\thanks{K. Yamamoto is with the Department of Electrical Engineering, Kyushu University, Fukuoka 819-0395, Japan.
}%
\thanks{This work was supported by the Swedish Research Council through the LCCC Linnaeus Center.}% <-this % stops a space
}

\begin{document}

\maketitle
\thispagestyle{empty}
\pagestyle{empty}

%%%%%%%%%%%%%%%%%%%%%%%%%%%%%%%%%%%%%%%%%%%%%%%%%%%%%%%%%%%%%%%%%%%%%%%%%%%%%%%%

\begin{abstract}
We present two sensitivity function trade-offs that apply to a class of networks with a string topology. In particular we show that a lower bound on the H-infinity norm and a Bode sensitivity relation hold for an entire family of sensitivity functions associated with growing the network. The trade-offs we identify are a direct consequence of growing the network, and can be used to explain why poorly regulated low frequency behaviours emerge in long vehicle platoons even when using dynamic feedback.
\end{abstract}

%%%%%%%%%%%%%%%%%%%%%%%%%%%%%%%%%%%%%%%%%%%%%%%%%%%%%%%%%%%%%%%%%%%%%%%%%%%%%%%%
\section{Introduction}

Feedback control is, at its heart, about managing trade-offs. Classically speaking, this is understood through the sensitivity function
\begin{equation}
S\s=\frac{1}{1+P\s{}C\s}.
\end{equation}
The `size' of the sensitivity function determines the effect of the feedback compensator $C\s$ on the open loop dynamics $P\s$. Trade-offs exist because $S\s$ cannot be freely assigned via the compensator. Instead, based on the model of the plant, the designer must choose which frequency ranges to improve, understanding that these improvements will come at the cost of degrading performance in other frequency ranges. 

Being able to manage such trade-offs is one of the major benefits of feedback control. This insight is central to many of the most successful model based synthesis methods, for example the loop-shaping method of McFarlane and Glover \cite{MG92}. However there is a fundamental shortcoming in these approaches when considering control problems in networks. The root of this issue is that in the majority of network applications, the designed controllers are required to work even as the network changes, for example as generators are added to a power grid, or as users join the Internet. Such changes affect the `plant model', and consequently the underlying trade-offs. Therefore it is necessary to explicitly consider the process of scaling the network model in order to properly design controllers for the network setting.

In order to formally address the above, we require a notion of scaling up a network. To this end, we define a sequence of plant and compensator models, and study the (1,1) entry of the sensitivity functions associated with this sequence
\[
S_N\s=\sqfunof{\funof{I+P_N\s{}C_N\s{}}^{-1}}_{1,1}.
\]
The meaning of this will be made precise in \Cref{sec:prob}, but the idea is that iterating this sequence corresponds to adding components to the network. Studying how this entry of the sensitivity function changes is then illustrative of how the effect of the feedback changes locally as the network grows. The precise setting we consider is highly simplified, but is suitable for networks with string topologies, such as vehicle platoons.

\begin{figure}[t]
\centering
\subfigure[$N=1$]{\input{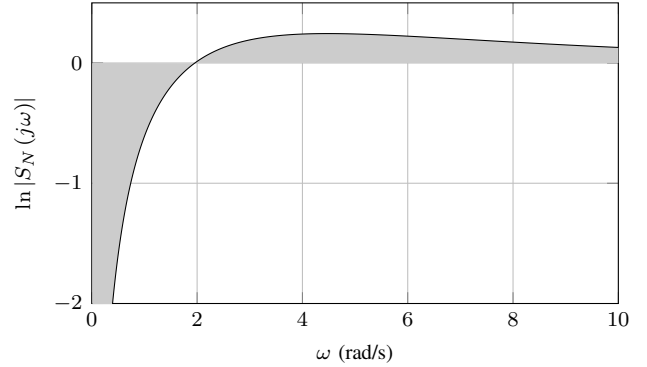}}\\
\subfigure[$N=10$]{\input{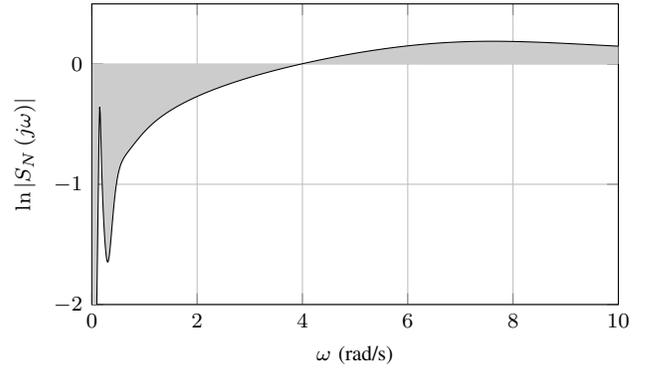}}\\
    \caption{Plots of $\abs{S_N\jw}$ for two values of $N$. Our first result is to show that peaks such as that near the origin in (b) can emerge as $N$ increases, and are an inevitable consequence of repeated imaginary axis poles in the plant. Our second is to show that a Bode-type integral relation applies to $S_N$, meaning that the (signed) shaded area in (a) and (b) is equal to zero, and remains invariant for all $N\in\mathbb{N}$.}
    \label{fig:macro}
\end{figure}

Given this notion of scaling a network, we consider network equivalents of two classical fundamental limitations on $S\s$. The phenomena we identify are directly related to the process of growing the network, and are illustrated in \Cref{fig:macro}.

\subsubsection{\Hfty{}-norm Limitations}

Tools such as Nevanlinna-Pick interpolation and the maximum modulus principle \cite{FZ84,dft92} can be used to derive lower bounds on the \Hfty{}-norm of $S\s$. That is they can be used to derive constants $k\geq{}0$ such that given any compensator,
\begin{equation}\label{eq:1}
\norm{S\s}_\infty\geq{}k.
\end{equation}
The values of the constants typically depend on the unstable poles and zeros of $P\s$, and impose fundamental limits on how small the sensitivity function can be made. Our first result is to show that the process of growing the network also imposes limits on the size of $S_N\s$. More specifically we show that there exists a constant $k^{\text{net}}$ such that
\[
\sup_{N\in\mathbb{N}}\norm{S_N\s}_\infty\geq{}k^{\text{net}}.
\]
The value of this constant depends on the residues of the unstable poles of the open loop plant, and imposes a fundamental limit on how small the sensitivity functions for the networks in the given class can be made. The restrictions are far more severe than in the scalar case, and can also result in the emergence of peaks in the sensitivity function near the open loop imaginary axis poles. This is precisely what we see in \Cref{fig:macro}(b), where the open loop plant has a repeated pole at the origin. This behaviour is in stark contrast to the scalar case, where $S\s$ is typically small near the open loop unstable poles. This, along with the connection to the large-scale accordion-like behaviours from \cite{BJM12,MB10}, is presented in \Cref{sec:hardy}.

\subsubsection{Bode Integral Relations}

Under mild assumptions, the Bode integral relation shows that there exists a constant $A\geq{}0$ that depends on the unstable poles of $P\s$ and $C\s$ such that
\begin{equation}\label{eq:2}
\int_0^\infty{}\ln\abs{S\jw}d\omega=A.
\end{equation}
This means that no compensator can uniformly reduce the size of $S\jw$, and making $S\jw$ smaller than 1 in some frequency range will necessarily make it larger in another. It is well known that similar relations hold for determinants of sensitivity functions in the multi-input multi-output setting \cite{Fre88}. We show that for the given network model, a similar relation holds for $S_N\s$. In particular, we demonstrate that if (to all intents and purposes) $S_N\s$ is stable for all $N$, then
\[
\funof{\forall{}N\in\mathbb{N}},\quad\int_{0}^\infty{}\ln\abs{S_N\jw}d\omega=0.
\]
That is a Bode-type integral relation is left invariant as the network grows, and the only possible value of the constant $A$ that will hold for all $N$ is zero. This means that the shaded regions in \Cref{fig:macro} have equal area. This is presented in \Cref{sec:bode}.

\section*{Notation}

The set of natural numbers is denoted by $\mathbb{N}$, and the open and closed right half plane by $\C_+$ and $\overline{\C}_+$ respectively. $\Rat$ denotes the set of real rational not necessarily proper transfer functions, and $\Hfty$ the Hardy space of transfer functions that are analytic on the open right half plane with norm $ \norm{G\s}_\infty \coloneqq \sup_{s\in\C_+}|G\s|$. The $(1,1)$ entry of a matrix $A$ is denoted by $[A]_{1,1}$.

\section{Problem Description}\label{sec:prob}

\begin{figure}[h]
\centering
\tikzstyle{block} = [draw, rectangle, minimum height=0.5cm, minimum width=1cm]
\tikzstyle{sum} = [draw, circle, node distance=1.5cm]
\tikzstyle{input} = [coordinate]
\tikzstyle{output} = [coordinate]
\tikzstyle{disturbance} = [coordinate]
\tikzstyle{pinstyle} = [pin edge={to-,thin,black}]
\begin{tikzpicture}[auto, node distance=1.5cm,>=latex']
    % We start by placing the blocks
    \node [input, name=input] {};
    \node [sum, right = 1cm of input] (sum1) {};
    \node [block, right of=sum1] (BT) {$B_N^T$};
    \node [block, right of=BT] (plant) {$P\s{}I$};
    \node [sum, right of=plant] (sum2) {};
    \node [output, right of=sum2] (output){};
    \node [disturbance, above = 1cm of sum2] (disturbance){};
    \node [block, below of=plant] (B){$B_N$};
    \node [block, below of=BT] (controller){$C\s{}I$};

    \draw [draw,->] (input) --  (sum1);
    \draw [->] (sum1) -- (BT);
    \draw [->] (BT) -- (plant);
    \draw [->] (plant) -- node[pos=0.99]{$+$}(sum2);
    \draw [->] (sum2) -- node[name=y]{$y$}(output);
    \draw [->] (disturbance)  --node{$d$}(sum2);
    \draw [->] (sum2)+(0.5,0) |- (B);
    \draw [->] (B) -- (controller);
    \draw [->] (controller) -| node[pos=0.99]{$-$}(sum1);
 \end{tikzpicture}
 \caption{The feedback configuration} 
 \label{fig:block}
\end{figure}
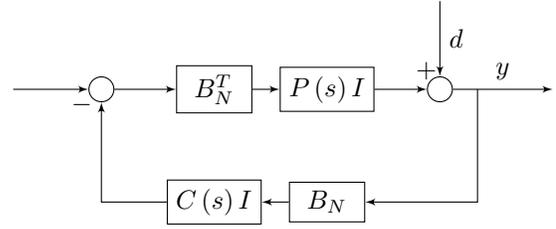
We consider the linear, time-invariant system depicted in \Cref{fig:block} where $P\s,C\s\in\Rat$ are the transfer functions of single-input single-output systems, $I$ is the identity matrix, and
\[
  B_N = \begin{bmatrix}
        1 & -1 & 0 &\cdots&0 \\
        0 & \ddots{} & \ddots{} &\ddots{} &\vdots{}\\
        \vdots{} & \ddots & \ddots & \ddots&0  \\
        \vdots{} & {} & \ddots{}& \ddots & -1\\
        0 & \cdots{} & \cdots & 0 &1
        \end{bmatrix} \in \R^{N\times N}.
\]
This feedback interconnection captures the dynamics of a network of $N$ homogeneous dynamical systems interconnected through a sparse homogeneous feedback law. This specific structure arises in a range of applications, including the symmetric bidirectional control of vehicle platoons \cite{SPH04,BHJ05}, and the control of a homogeneous chain of one-dimensional harmonic oscillators \cite{YS16}. It also applies to frequency and voltage stability problems in electrical power systems \cite{MBB97}, and flocking and consensus phenomena \cite{BJM12}, under the restriction that the network has a line structure.

Our particular focus is to study how the effect of the feedback changes locally as the network grows. To this end we introduce the following infinite family of transfer functions
\[
S_N\s\coloneqq{}\sqfunof{(I+P\s{}C\s{}L_N)^{-1}}_{1,1},
\]
where
\[
    L_N \coloneqq B_N^TB_N =
         \begin{bmatrix}
        1 & -1 & 0 &\cdots&0 \\
        -1 & 2 & \ddots{} &\ddots{} &\vdots{}  \\
        0 & \ddots & \ddots & \ddots&0  \\
        \vdots{}&\ddots{}&\ddots{}&\ddots{}&-1\\
        0 &  \cdots{}    &0& -1 & 2
        \end{bmatrix} \in \R^{N\times N}.
\]
This is nothing but the \funof{1,1} entry of the output sensitivity function (the transfer function from $d$ to $y$) in a system of size $N$. This transfer function therefore captures how a local disturbance affects the local behaviour of the system; more specifically how $d_1$ affects $y_1$. A method to design the compensator $C(s)$ in a scale free manner with respect to $S_N\s$ has been presented in \cite{PY18}. In the sequel we study fundamental limitations on this design that are imposed by the properties of $P\s$.

\section{Results}

\subsection{\Hfty{}-norm Limitations}\label{sec:hardy}

In this subsection we will present a lower bound on $\norm{S_N\s}_\infty$. This lower bound is obtained by considering the behaviour of $S_N\s$ around the open loop poles of $P\s{}C\s$, and has two main consequences:
\begin{enumerate}
    \item If $P\s$ contains an unstable pole in the open right-half plane, then no compensator can bound the \Hfty{}-norm of $S_N\s$ for all $N\in\mathbb{N}$.
    \item If $P\s$ contains an imaginary axis pole of multiplicity 2, then there is a lower bound on $\norm{S_N\s}_\infty$ that depends on the residues of the pole.
\end{enumerate}
At first sight, point 2 appears highly specialised. However, this is precisely the setting considered in the vehicle platooning literature, in which the plant has a repeated pole at the origin (see e.g. \cite{SPH04,BHJ05}). This lower bound therefore imposes a restriction on the effect of any compensator's ability to regulate low frequency behaviours as the network size becomes large. This manifests itself through a peak in the sensitivity function around the open loop pole, and is highly suggestive as to the fundamental nature of the accordion type behaviours of \cite{BJM12} even when dynamical controllers are used. 

\begin{figure}
\centering
    \input{Figures/SensitivityBode.tex}
    \caption{Plot of $\abs{S_N\jw}$ for different values of $N$, when $P\s{}C\s$ is given by \cref{eq:ex2}.}
    \label{fig:2}
\end{figure}

\begin{theorem}\label{thm:hardy}
Let $P\s,C\s\in\Rat$, and assume that $P\s{}C\s$ has a pole of multiplicity $m$ at $p\in\overline{\C}_+$. Define the Laurent series expansion of ${P\s{}C\s}$ at $s=p$ as
\begin{equation}\label{eq:laurent}
 P\s{}C\s=\sum_{k=-\infty}^\infty{}a_k\funof{s-p}^k.
\end{equation}
If $p\in{}j\R$ and $m\leq{}2$, then
\[
\sup_{N\in\mathbb{N}}\norm{S_N\s}_\infty\geq{}\frac{4\funof{m-1}}{\pi{}\abs{{a_{1-m}\sqrt{a_{-m}}}}},
\]
otherwise
\[
\sup_{N\in\mathbb{N}}\norm{S_N\s}_\infty=\infty.
\]
\end{theorem}

Before proving the result we give two examples of its application. The first explains the emergence of the peak in \Cref{fig:macro} when a vehicle platoon model is used, and the second shows the same behaviour appearing when $P\s{}C\s$ contains a repeated pole at $s=j$.

\begin{example}\label{ex:1}
Suppose that
\[
P\s{}=\frac{1}{s^2\funof{0.1s+1}},C\s=\frac{2s+1}{0.05s+1}.
\]
This is the plant and controller used to illustrate the results on vehicle platoons from \cite{SPH04}, and is also the model we used to generate \Cref{fig:macro}. 
This transfer function contains a pole of multiplicity 2 at the origin, with Laurent series expansion given by
\[
P\s{}C\s=\frac{1}{s^2}+\frac{7}{4s}+\cdots{}.
\]
Therefore by \Cref{thm:hardy},
\[
\sup_{N\in\mathbb{N}}\norm{S_N\s}_\infty{}\geq{}\frac{16}{7\pi{}}.
\]
The emergence of the peak around the origin is illustrated in \Cref{fig:macro}. Note that the height of the peak is $\ln\funof{\nicefrac{16}{7\pi{}}}\approx{}-0.3$, as predicted by \Cref{thm:hardy}. This shows that the compensator does not effectively regulate low frequency disturbances. We note the similarity of this observation to that involving the accordion like behaviour in platoons from \cite{BJM12} and the unavoidable long transient period under the bidirectional control scheme from \cite{MB10}. There they demonstrated that under simple control strategies a poorly regulated low frequency behaviour emerges as the number of vehicles in a platoon becomes large. Here we see precisely the same thing happening even when more complex dynamical controllers are considered.
\end{example}
\begin{example}\label{ex:2}
Suppose that
\begin{equation}\label{eq:ex2}
P\s{}C\s=\frac{\funof{s+1}^4}{\funof{s^2+1}^2}.
\end{equation}
$P\s{}C\s$ has two poles of multiplicity 2 at $s=\pm{}j$. The Laurent series expansion of $P\s{}C\s$ at $s=j$ is given by
\[
P\s{}C\s{}=\frac{1}{\funof{s-j}^2}+\frac{2-j}{s-j}+\ldots{}
\]
Therefore by \Cref{thm:hardy},
\[
\sup_{N\in\mathbb{N}}\norm{S_N\s}_\infty{}\geq{}\frac{4}{\pi{}\sqrt{5}}\approx{}-5 \,\mathrm{dB}.
\]
This is illustrated in \Cref{fig:2}. As $N$ increases we see the emergence of a peak in $S_N\jw$ around the open loop pole with a height of $\approx{}-5\,\mathrm{dB}$.
\end{example}

\begin{proof}
The proof will be in two main stages. First we will show that if $P\s$ has any open right half plane poles, or imaginary axis poles of multiplicity three or greater, then there exists no controller such that $S_N\s\in\Hfty{}$ for all $N$. We will then show that when the imaginary axis poles have multiplicity two or less, the bound from the first part of the theorem holds. 

\textit{Stage 1:}
The poles of $S_N\s$ are given by the solutions to the equation
\[
\det\funof{I+P\s{}C\s{}L_N}=0.
\]
Therefore if we can show that for some eigenvalue $\lambda$ of $L_N$ and some $u\in\overline{\C}_+$
\[
-\frac{1}{P\funof{u}C\funof{u}}=\lambda,
\]
then $S_N\s$ will have a pole at the point $u$, and hence be unstable. As $N$ becomes large, the eigenvalues of $L_N$ become dense on the interval \funof{0,4} (see \cref{eq:eigeig} in the proof of \Cref{thm:bode} for an expression). Consequently a sufficient condition for instability is that there exists a set $U\subset\overline{\C}_+$ and numbers $0<\epsilon_1<\epsilon_2<4$ such that
\[
\sqfunof{\epsilon_1,\epsilon_2}\in\cfunof{v:v=-\frac{1}{P\funof{u}C\funof{u}},u\in{}U}.
\]
Consider now the Laurent series expansion in \cref{eq:laurent}, which shows that
\begin{equation}\label{eq:laurent1}
\frac{1}{P\s{}C\s}=\frac{\Delta{}s^m}{a_{-m}}+\mathcal{O}\funof{\Delta{}s^{m+1}}.
\end{equation}
In the above $\Delta{}s=\funof{s-p}$. By considering the annulus
\[
U=\cfunof{u:r_1<\abs{u-p}<r_2},
\]
if $0<r_1<r_2$ are sufficiently small, it can be seen directly from \cref{eq:laurent1} that this annulus is mapped by $\frac{1}{P\s{}C\s}$ into a region that is arbitrarily close to the origin, and encircles it. Therefore this new region will contain a suitable interval $\sqfunof{\epsilon_1,\epsilon_2}$, and consequently $S_N\s$ will contain a pole in $\overline{\C}_+$ for some $N$. A similar argument can be used in the case of imaginary axis poles. This time we cannot use an annulus for the set $U$, since this will not lie in $\overline{\C}_+$. However, we can use a set of the form
\[
U=\cfunof{u:r_1<\abs{u-p}<r_2,\text{Re}\funof{u}>0}.
\]
This time provided $m>2$, $U$ will be mapped into a region that encircles the origin, again showing that $S_N\s$ is unstable for some $N$.

\textit{Stage 2:} We require the following expression for $S_N\s$, which comes from the theory of iterated M\"{o}bius transformations, (this specific equation is from \cite[\S{}V.B]{YS16}):
\begin{equation}\label{eq:multform}
\begin{aligned}
\Sn\s&=\funof{1-\zeta\s}\frac{1-\zeta\s^{2{N}}}{1+\zeta\s^{2N+1}}.
\end{aligned}
\end{equation}
In the above $\zeta\s$ is given by the root of
\begin{equation}\label{eq:mult}
\zeta\s^2-\funof{\frac{1}{P\s{}C\s}+2}\zeta\s+1=0
\end{equation}
satisfying $\abs{\zeta\s}<1$. The bound we derive comes from the behaviour of $S_N\s$ in the neighbourhood of the pole $p$, so once again we work with the Laurent series expansion. From \cref{eq:laurent} it follows that along the imaginary axis
\[
\frac{1}{P\jw{}C\jw}=-\frac{\Delta{}\omega^2}{a_{-2}}+\frac{ja_{-1}\Delta{}\omega^3}{a_{-2}}+\mathcal{O}\funof{\Delta{}\omega^4},
\]
where $\Delta\omega=\funof{j\omega-p}.$ It can then be shown using the Maclaurin series for $\sqrt{1+z}$ that the solution to \cref{eq:mult} satisfies
\[
\zeta\jw=1-\frac{j\Delta{}\omega}{\sqrt{a_{-2}}}+\frac{\funof{a_{-1}\sqrt{a_{-2}}-1}\Delta{}\omega^2}{2a_{-2}}+\mathcal{O}\funof{\Delta{}\omega^3}.
\]
Now consider the sequence of frequencies
\[
\Delta\omega_N=\frac{\pi\sqrt{a_{-2}}}{{2N+1}}.
\]
It then follows that
\[
\begin{aligned}
1-\zeta\funof{j\omega_N}&=\frac{j\pi}{2N+1}+\mathcal{O}\funof{\Delta{}\omega_N^2},\\
1-\zeta\funof{j\omega_N}^{2N}&=2+\mathcal{O}\funof{\Delta{}\omega_N},\\
1+\zeta\funof{j\omega_N}^{2N+1}&=\frac{-\pi^2}{2N+1}\funof{\frac{a_{-1}\sqrt{a_{-2}}}{2}}+\mathcal{O}\funof{\Delta{}\omega_N^2}.
\end{aligned}
\]
This shows that
\[
\lim_{N\rightarrow\infty}S_N\funof{j\omega_N}=\frac{-4j}{\pi{}{a_{-1}\sqrt{a_{-2}}}}.
\]
This implies the bound for $m=2$, because
\[
\sup_{N\in\mathbb{N}}\norm{S_N}_\infty\geq{}\sup_{N\in\mathbb{N}}\abs{S_N\funof{j\omega_N}}.
\]
The bound for $m=1$ is trivial (it can also be shown that running the above argument with $m=1$ results in the trivial lower bound of zero).
\end{proof}

\subsection{Sensitivity Integral Relations}\label{sec:bode}

In this subsection we present a Bode-type integral relation for the family of transfer functions $S_N\s$. The result shows that if a test involving $P\s{}C\s$ is satisfied, then a sensitivity integral relation is left invariant for all $N\in\mathbb{N}$. This test can be conducted with standard tools such as the Routh-Hurwitz stability criterion, and is to all intents and purposes equivalent to $S_N\s$ being stable for all $N\in\mathbb{N}$.

\begin{theorem}\label{thm:bode}
Let $P\s,C\s\in\Rat$. If
\[
\funof{\forall{}k\in\funof{0,4}},\quad\frac{1}{1+kP\s{}C\s}\in\Hfty{},
\]
then $S_N\s\in\Hfty$ for all $N\in\mathbb{N}$. If in addition $P\s{}C\s$ has at least two more poles than zeros, then
\[
\funof{\forall{}N\in\mathbb{N}}, \quad \int_0^\infty \ln |S_N(j\omega)| d\omega= 0.
\]
\end{theorem}

Before proving the result we give an example of its application.

\begin{example}\label{ex:3}
Consider again $P\s{}C\s$ from \Cref{ex:1}. 
The denominator of $\funof{1+kP\s{}C\s}^{-1}$ is
\begin{equation}\label{eq:poly}
k\funof{2s+1}+s^2\funof{0.1s+1}\funof{0.05s+1}.
\end{equation}
By the Routh-Hurwitz stability criterion, a fourth order polynomial $a_4s^4+\ldots{}+a_0$ is stable if and only if $a_i>0$ and $a_1a_2a_3>a_0a_3^2+a_4a_1^2$. Consequently \cref{eq:poly} is stable for all $13\,\nicefrac{7}{8}>k>0$, and hence $S_N\s$ is stable for all $N\in\mathbb{N}$ by \Cref{thm:bode}. Furthermore, since $P\s{}C\s$ has 1 zero and 4 poles, the integral relation in \Cref{thm:bode} also holds. This is illustrated in \Cref{fig:macro}.
\end{example}

\begin{proof}
We first prove the claim about the stability of the functions $S_N\s$, and then show that the integral relation holds.

\textit{Step 1:} $S_N\s$ is stable if and only if
\[
\begin{aligned}
\frac{1}{\det\funof{I+P\s{}C\s{}L_N}}\in\Hfty.
\end{aligned}
\]
The eigenvalues of $L_N$ are given by 
\begin{equation}\label{eq:eigeig}
    2\funof{1-\cos\frac{(2k-1)\pi}{2N+1}}, \quad k=1,2,\dots,N.
\end{equation}
Since the determinant of a matrix is equal to the product of its eigenvalues, $S_N\s$ is stable if and only if
\begin{equation}\label{eq:fin}
{\prod_{k=1}^N\frac{1}{\funof{1+2\funof{1-\cos\frac{(2k-1)\pi}{2N+1}}P\s{}C\s}}}\in\Hfty.
\end{equation}
Since $2\funof{1-\cos\frac{(2k-1)\pi}{2N+1}}\in(0,4)$, the conditions of the theorem imply \cref{eq:fin}, and consequently that $S_N\s\in\Hfty$ for all $N\in\mathbb{N}$ as required.

\textit{Step 2:} Now define 
\[
\overline{L}_N = \begin{bmatrix}
        2 & -1 & 0 &\cdots&0 \\
        -1 & 2 & \ddots{} &\ddots{} &\vdots{}  \\
        0 & \ddots & \ddots & \ddots&0  \\
        \vdots{}&\ddots{}&\ddots{}&\ddots{}&-1\\
        0 &  \cdots{}    &0& -1 & 2
        \end{bmatrix} \in \R^{N\times N}.
\]
Observe that the (1,1) element of $\mathrm{adj}(I+P\s{}C\s{}L_N)$ is 
${\det(I+P\s{}C\s\overline{L}_{N-1}).}$
Hence (supressing the dependence on $s$), 
\[
\begin{aligned}
\sqfunof{(I+PCL_N)^{-1}}_{1,1} &= 
\frac{\sqfunof{\mathrm{adj}\funof{I+PCL_N}}_{1,1}}{\det\funof{I+PCL_N}} \\
        &= \frac{\det(I+PC\overline{L}_{N-1})}{\det(I+PCL_N)}.
\end{aligned}
\]
Therefore, 
\[
\begin{aligned}
    \int_0^\infty \ln |S_N(j\omega)| d\omega
    ={}& \int_0^\infty \ln \left|\frac{\det(I+PC\overline{L}_{N-1})}{\det(I+PCL_N)}\right|d\omega \\
    ={}& \int_0^\infty \ln |\det(I+PCL_N)^{-1}|d\omega \\
    &- \int_0^\infty \ln |\det(I+PC\overline{L}_{N-1})^{-1}|d\omega\\
\end{aligned}
\]
This final equation is of the correct form to be evaluated using the sensitivity integral relation for multivariable systems given in \cite[Theorem 5.3.1]{Fre88}. This result states that given any transfer matrix $M\s\in\Rat^{N\times{}N}$ with $N_p$ poles at $p_i \in \overline{\C}_+$, $i=1,\dots N_p$, if $\funof{I+M\s}^{-1}\in\Hfty^{N\times{}N}$, then
\[
\int_0^\infty{}\ln\abs{\det\funof{I+M(j\omega)}^{-1}}d\omega=\pi \sum_{i=1}^{N_p} \mathrm{Re}\funof{p_i}.
\]
We proved in Step 1 that $\funof{I+P\s{}C\s{}L_N}^{-1}$ is stable for all $N\in\mathbb{N}$. In fact, since the eigenvalues of $\overline{L}_N$ are given by 
\[
    2\funof{1-\cos\frac{k\pi}{N+1}}, \quad k=1,2,\dots,N,
\]
an identical argument shows that $\funof{I+P\s{}C\s{}\bar{L}_N}^{-1}$ is stable for all $N\in\mathbb{N}$. Therefore
\[
\int_0^\infty \ln |S_N(j\omega)| d\omega =\pi \sum_{i=1}^{N_p} \mathrm{Re}\funof{p_i},
\]
where $p_i$ denote the unstable poles of $P\s{}C\s$. Finally we note that by \Cref{thm:hardy}, if $S_N\s$ is stable for all $N$ then $P\s{}C\s$ can have no poles in the open right half plane, and hence $\sum_{i=1}^{N_p}\text{Re}\funof{p_i}=0$. This completes the proof.
\end{proof}

\section{Conclusions}
Two sensitivity function trade-offs that apply to networks with homogeneous agent dynamics and a string topology have been presented. The first shows that the residues of the open loop unstable poles of the agent dynamics impose a lower bound on the $\Hfty{}$-norm of the sensitivity function. In particular, an open-loop imaginary axis pole of multiplicity $2$ results in the emergence of a peak in the sensitivity function as the network grows. Moreover, no controller can bound the $\Hfty{}$-norm of the sensitivity function independently of network size if an open-loop unstable pole in the open right-half plane exists. It has also been shown that a Bode-type integral relation holds for the studied sensitivity function. This means that the `waterbed effect' is present in networks of any size in a manner entirely analogous to the single-input single-output case. Several numerical examples have been given to illustrate the results.

\bibliographystyle{IEEEtran}
\bibliography{references.bib}

\end{document}

%% file: Maths.tex
\newcommand{\abs}[1]{\ensuremath{\left\vert#1\right\vert}}

\newcommand{\cfunof}[1]{\ensuremath{\left\{#1\right\}}}

\newcommand{\C}{\ensuremath{\mathbb{C}}}

\newcommand{\funof}[1]{\ensuremath{\left(#1\right)}}

\newcommand{\Hfty}{\ensuremath{\mathscr{H}_{\infty}}}

\newcommand{\jw}{\ensuremath{\left(j\omega\right)}}

\newcommand{\norm}[1]{\ensuremath{\left\Vert #1 \right\Vert}}

\newcommand{\Rat}{\ensuremath{\mathscr{R}}}

\newcommand{\R}{\ensuremath{\mathbb{R}}}

\newcommand{\s}{\ensuremath{\left(s\right)}}
\newcommand{\sqfunof}[1]{\ensuremath{\left[#1\right]}}

\newcommand{\Sn}{\ensuremath{S_N}}